\documentclass[11pt,]{article}
\usepackage[latin1]{inputenc}
\usepackage[english]{babel}
\usepackage{amsfonts}
\usepackage{amssymb}
\usepackage{amsmath}
\usepackage{latexsym}
\usepackage{amsthm}
\usepackage[initials]{amsrefs}
\newtheorem{theorem}{Theorem}[section]

\newtheorem{proposition}[theorem]{Proposition}
\newtheorem{lemma}[theorem]{Lemma}
\theoremstyle{definition}
\newtheorem{definition}[theorem]{Definition}

\newtheorem{remark}[theorem]{Remark}
\newenvironment{keyw}{\textit{Keywords:}}

\newenvironment{MSC}{Mathematics Subject Classification 2000:} 

\newcommand{\cy}[2]{C(#1,#2)}
\newcommand{\cyl}[2]{C_{#1}({#2})}
\DeclareMathOperator{\aut}{Aut}
\DeclareMathOperator{\spa}{span}
\DeclareMathOperator{\Id}{Id}

\def\CC{\mathcal{C}}
\def\C{\mathbb{C}}
\def\B{\mathsf{L}}
\def\A{\mathcal{A}}
\def\E{\mathcal{E}}
\def\D{\mathcal{D}}
\def\Dt{\widetilde{\mathcal{D}}}
\def\Oo{\mathcal{O}}
\def\L{\mathcal{L}}

\def\K{\mathcal{K}}

\def\J{\mathcal{J}}

\def\N{\mathbb{N}}
\newcommand{\No}{{\N_0}}
\def\Z{\mathbb{Z}}
\def\T{\mathbb{T}}
\def\a{\mathfrak{a}}

\def\X{\mathcal{X}}

\def\H{\mathsf{H}}

\newcommand{\kd}{k_{\Dt_{\OSS}}}
\newcommand{\kh}{k_{\H_{\OSS}}}
\newcommand{\emptyword}{\epsilon}

\newcommand{\cs}{C^*}

\newcommand{\OSS}{{\mathsf{X}}}
\newcommand{\OSSY}{{\mathsf{Y}}}

\newcommand{\osh}{\sigma}
\newcommand{\tsh}{\sigma}
\newcommand{\inn}[1]{\langle {#1}\rangle}
\newcommand{\spac}{\overline{\spa}}
\providecommand{\norm}[1]{\lVert#1\rVert}
\author{Toke Meier Carlsen\footnote{Current address: Mathematisches
    Institut, Einsteinstra\ss e 62, 48149 M\"unster, Germany}}
\title{Cuntz-Pimsner $C^*$-algebras associated with subshifts}
\date{}

\begin{document}

\maketitle
\begin{center}
  \textit{Institut for Matematiske Fag, K\o benhavns Universitet}\\
  \textit{Universitetsparken 5, 2100 K\o benhavn \O{}, Denmark}\\
  \textit{toke@math.ku.dk}
  
  \vspace{5mm}
  
\end{center}

\begin{abstract}
  \noindent By using $\cs$-correspondences and Cuntz-Pimsner algebras, we
  associate to every subshift (also called a shift space) $\OSS$ a
  $\cs$-algebra $\Oo_{\OSS}$, which is a generalization of the
  Cuntz-Krieger algebras. We show that $\Oo_{\OSS}$ is the universal
  $\cs$-algebra generated by partial isometries satisfying relations
  given by $\OSS$. We also show that $\Oo_{\OSS}$ is a one-sided conjugacy
  invariant of $\OSS$.
  \vspace{4mm}

  \noindent
  \begin{keyw}
    $C^*$-algebras, subshifts, shift spaces, conjugacy, Cuntz-Krieger
    algebras, Cuntz-Pimsner algebras. 
  \end{keyw}
  \vspace{4mm}

  \noindent
  \begin{MSC}
    Primary: 46L55, Secondary: 37B10.
  \end{MSC}
\end{abstract}

\section{Introduction}
In \cite{MR82f:46073a} Cuntz and Krieger introduced a new class of
$C^*$-algebras which in a natural way can be viewed as universal $C^*$-algebras
associated with subshifts (also called shift spaces) of finite
type. From the point of view of
operator algebra these $C^*$-algebras were important examples of
$C^*$-algebras with new properties and from the point of view of
topological dynamics these $C^*$-algebras (or rather, the $K$-theory
of these $C^*$-algebras) gave new invariants of subshifts of
finite type.

In \cite{MR98h:46077} Matsumoto tried to generalize this idea by constructing
$C^*$-algebras associated with every subshift, and he and others have
studied these $C^*$-algebras in
\cites{MR1978330,MR2000e:46087,MR2000d:46082,MR2000f:46084,MR2001e:46115,MR2001g:46147,MR1645334,MR1637788,MR1716953,MR2002h:19004}.
Unfortunately there is a mistake in \cite{MR2000f:46084} which makes
many of the results in
\cites{MR2000d:46082,MR2000f:46084,MR2001e:46115,MR2001g:46147,MR1637788}
invalied for the $C^*$-algebras constructed in
\cite{MR98h:46077}. This mistake has to do with the identification 
of an underlying compact space which among other things determine the
$K$-theory of the $C^*$-algebras. It turned out that this compact
space is not the space Matsumoto thought it was, and thus many of the
results of
\cites{MR2000d:46082,MR2000f:46084,MR2001e:46115,MR2001g:46147,MR1637788}
are invalied for the $C^*$-algebras constructed in
\cite{MR98h:46077}. To recover these results Matsumoto and the
author introduced in \cite{MR2091486} a new class of $C^*$-algebras
associated with subshifts, which has the right underlying compact
space and thus satisfies most of the results in
\cites{MR1978330,MR2000e:46087,MR98h:46077,MR2000d:46082,MR2000f:46084,MR2001e:46115,MR2001g:46147,MR1645334,MR1637788,MR1716953,MR2002h:19004},
but these $C^*$-algebras do not in general
have the universal property (cf. Theorem \ref{unihomo}, Remark
\ref{remarket} and \ref{tocalg}, \cite{MR98h:46077}*{Theorem 4.9} and
\cite{MR2091486}*{pp. 148-149}). Thus it is natural to think of this
class of $\cs$-algebras as the class of 
reduced $\cs$-algebras associated with subshifts.

In this paper we will for each subshift $\OSS$ construct a
new $C^*$-algebra $\Oo_{\OSS}$ by using $\cs$-correspon\-dences (also
called Hilbert 
bimodules) and Cuntz-Pimsner algebras, and this new $C^*$-algebra
will both have the right underlying compact space and have the
universal property and hence will satisfy all the results of
\cites{MR1978330,MR2000e:46087,MR98h:46077,MR2000d:46082,MR2000f:46084,MR2001e:46115,MR2001g:46147,MR1645334,MR1637788,MR1716953,MR2002h:19004}
and has the $\cs$-algebra defined in \cite{MR2091486} as a quotient.
Thus is seems right to think of this $\cs$-algebra as the universal
$\cs$-algebra associated to a subshift. 

The $\cs$-algebra $\Oo_{\OSS}$ can also be constructed as the
$\cs$-algebra of a groupoid (cf. Remark \ref{remark:groupoid}), and by
using Exel's crossed product of a $\cs$-algbra of an endomorphism (cf. Remark
\ref{remark:sergei}). 

Matsumoto's original construction associated a $\cs$-algebra to every
\emph{two-sided} subshift, but it seems more natural to work with
\emph{one-sided} subshifts, so we will do that in this paper, but
since every two-sided subshift comes with a canonical one-sided
subshift (see below), the $\cs$-algebras we define in this paper can
in a natural way also been seen as $\cs$-algebras associated to
two-sided subshifts (cf. Remark \ref{tocalg}). 

We will show that $\Oo_{\OSS}$ is the universal $\cs$-algebra
associated with partial isometries satisfying relations giving by
$\OSS$ and which resemble the Cuntz-Krieger relations (Theorem
\ref{unihomo}). We will also show that $\Oo_{\OSS}$ is an invariant
of $\OSS$ in the sense that if $\OSS$ and $\OSSY$ are conjugate
one-sided subshifts, then $\Oo_{\OSS}$ and $\Oo_{\OSSY}$ are
isomorphic (Theorem \ref{toke}). This is a generalization of
\cite{MR82f:46073a}*{Proposition 2.17} and 
\cite{MR98h:46077}*{Proposition 5.8} (see \cite{MR2000d:46082}*{Lemma
  4.5} for a proof of the later Proposition), where it is required
that $\OSS$ and $\OSSY$ satisfy a certain condition (I). 

\section{Notation and preliminaries} \label{notation}
Throughout this paper, $\No$ will denote the set of non-negative
integers.

Let $\a$ be a finite set endowed with the discrete topology. We will
call this set the alphabet. Let $\a^\No$ be the infinite
product spaces $\prod_{n=0}^\infty \a$
endowed with the product topology. The transformation
$\tsh$ on $\a^\No$ given by $(\tsh(x))_i=x_{i+1},\ i\in
\No$ is called the \emph{shift}. Let $\OSS$ be a shift invariant
closed subset of $\a^\No$ (by shift invariant we mean that
$\tsh(\OSS)\subseteq \OSS$, not necessarily
$\tsh(\OSS)=\OSS$). The topological dynamical system
$(\OSS,\tsh_{|\OSS})$ is called a \emph{subshift}. We will denote
$\tsh_{|\OSS}$ by $\tsh_{\OSS}$ or $\tsh$ for simplicity,
and on occasion the alphabet $\a$ by $\a_{\OSS}$. Since $\tsh_{\OSS}$
maps $\OSS$ into $\OSS$, we can compose $\tsh_{\OSS}$ with itself. We
will denote this map by $\tsh^2_{\OSS}$ and in general for any
positive integer $n$ the $n$-fold
composition of $\tsh_{\OSS}$ with itself by $\tsh^n_{\OSS}$. We will
for a subset $Y$ of the subshift $\OSS$ and an integer $n$ by
$\osh^n(Y)$ denote 
\begin{equation*}
  \osh^n(Y)=
  \begin{cases}
    \osh^n_{\OSS}(Y)& n>0,\\
    Y& n=0,\\
    (\osh^{-n}_{\OSS})^{-1}(Y)& n<0.
  \end{cases}
\end{equation*}


A finite sequence $\mu=(\mu_1,\ldots ,\mu_k)$ of elements $\mu_i\in
\a$ is called a finite \emph{word}. The \emph{length} of $\mu$ is $k$ and is denoted 
by $|\mu|$. We let for each $k\in \No$, $\a^k$ be the set of all words
with length $k$ and we let
$\B^k(\OSS)$ be the set of all words with length $k$ appearing in
some $x\in \OSS$. We set $\B_l(\OSS)=\bigcup_{k=0}^l
\B^k(\OSS)$ and $\B(\OSS)=\bigcup_{k=0}^\infty \B^k(\OSS)$ and likewise
$\a_l=\bigcup_{k=0}^l\a^k$ and $\a^*=\bigcup_{k=0}^\infty \a^k$, 
where $\B^0(\OSS)=\a^0$ denote the set consisting of the empty word
$\emptyword$ which has length 0. $\B(\OSS)$ is called the \emph{language} of $\OSS$. Note
that $\B(\OSS)\subseteq \a^*$ for every subshift.

For a subshift $\OSS$ and a word $\mu \in \B(\OSS)$ we denote by $\cyl{\OSS}{\mu}$ the \emph{cylinder set} 
\begin{equation*}
\cyl{\OSS}{\mu}=\{x\in \OSS \mid (x_1,x_2,\dots,x_{|\mu|})=\mu\}.
\end{equation*}
It is easy to see that
\begin{equation*}
\{\cyl{\OSS}{\mu}\mid \mu \in \B(\OSS)\}
\end{equation*}
is a basis for the topology of $\OSS$, and that $\cyl{\OSS}{\mu}$ is
closed and hence compact for every $\mu\in \B(\OSS)$.
We will allow us self to write $C(\mu)$ instead of $\cyl{\OSS}{\mu}$
when it is clear which subshift space we are working with.

For a subshift $\OSS$ and words $\mu, \nu \in \B(\OSS)$ we
denote by $\cy\mu\nu$ the set 
\begin{equation*}
C(\nu)\cap\tsh^{-|\nu|}(\tsh^{|\mu|}(C(\mu)))= \{\nu x\in \OSS\mid \mu
x\in \OSS\}.
\end{equation*}

If $\OSS$ and $\OSSY$ are two subshifts and $\phi:\OSS\to\OSSY$ is a
homeomorphism such that
$\psi\circ{\osh}_{\OSS}={\osh}_{\OSSY}\circ\phi$, then we say that $\phi$
is a \emph{conjugacy} and that $\OSS$ and $\OSSY$ are \emph{conjugate}.

What we have defined above is a \emph{one-sided} subshift. A
\emph{two-sided} subshift is defined in the same way, except that we
replace $\No$ with $\Z$: Let $\a^\Z$ be the infinite
product spaces $\prod_{n=-\infty}^\infty \a$
endowed with the product topology, and let $\tsh$ be the transformation
on $\a^\Z$ given by $(\tsh(x))_i=x_{i+1},\ i\in
\Z$. A shift invariant closed subset
$\Lambda$ of $\a^\Z$  (here, by shift
invariant we mean $\tsh(\Lambda)=\Lambda$) is called a \emph{two-sided
  subshift}. The set 
\begin{equation*}
  \OSS_\Lambda=\{(x_i)_{i\in\No}\mid (x_i)_{i\in\Z}\in \Lambda\}
\end{equation*}
is a one-sided subshift, and it is called the one-sided subshift of
$\Lambda$.

\section{Cuntz-Pimsner algebras}
We will in this section give a short introduction to Cuntz-Pimsner
algebras. We will follow the universal approach of \cite{MR1986889}
(see also \cite{MR97k:46069}, \cite{MR2002f:46113} and \cite{katsura}).

  Let $\A$ be a $C^*$-algebra. A right \emph{Hilbert $\A$-module} $\H$ is a
  Banach space with a right action of the $\cs$-algebra $\A$ and an
  $\A$-valued inner product $\inn{\cdot,\cdot}$ satisfying
  \begin{enumerate}
  \item $\inn{\xi,\eta a}=\inn{\xi,\eta}a,$
  \item $\inn{\xi,\eta}=\inn{\eta,\xi}^*$,
  \item $\inn{\xi,\xi}\ge 0$ and $\norm{\xi}=\norm{\inn{\xi,\xi}}^{1/2}$,
  \end{enumerate}
  for $\xi,\eta\in \H$ and $a\in \A$.

  For a Hilbert $\A$-module $\H$, we denote by $\L(\H)$ the $\cs$-algebra
  of all adjointable operators on $\H$. For $\xi,\eta\in \H$, the
  operator $\theta_{\xi,\eta}\in \L(\H)$ is defined by
  $\theta_{\xi,\eta}(\zeta)=\xi\inn{\eta,\zeta}$ for $\zeta\in \H$. We
  define $\K(\H)\subseteq \L(\H)$ by
  \begin{equation*}
    \K(\H)=\spac\{\theta_{\xi,\eta}\mid \xi,\eta\in \H\},
  \end{equation*}
  where $\spac\{\cdots\}$ means the closure of the linear span of $\{\cdots\}$.

Let $\phi:\A\to \L(\H)$ be a $*$-homomorphism. Then $ax:=\phi(a)x$ defines a left
action of $\A$ on $\H$, and we call $\H$ a \emph{$\cs$-correspondence} over $\A$
(in \cite{MR97k:46069} and \cite{MR1986889} a $\cs$-correspondence is called a
Hilbert bimodule, but it now seems that the term Hilbert bimodule has been
reserved for a special kind of $\cs$-correspondences cf. \cite{MR2002f:46113}). 

A \emph{Toeplitz representation} $(\psi,\pi)$ of $\H$ in a $C^*$-algebra $B$
consists of a linear map $\psi:\H\to B$ and a $*$-homomorphism
$\pi:\A\to B$ such that 
\begin{equation*}
  \psi(\xi a)=\psi(\xi)\pi(a),\ \psi(\xi)^*\psi(\eta)=\pi(\inn{\xi,\eta}),\text{
    and } \psi(a\xi)=\pi(a)\psi(\xi)
\end{equation*}
for $\xi,\eta\in \H$ and $a\in \A$. Given such a representation, there is a
homomorphism $\pi^{(1)}:\K(\H)\to B$ which satisfies
\begin{equation*}
  \pi^{(1)}(\theta_{\xi,\eta})=\psi(\xi)\psi(\eta)^*
\end{equation*}
for all $\xi,\eta\in \H$, and we then
have 
\begin{equation*}
  \pi^{(1)}(T)\psi(\xi)=\psi(T\xi)
\end{equation*}
for every $T\in \K(\H)$ and $\xi\in \H$. If $\rho:B\to C$ is a
$*$-homomorphism between $C^*$-algebras, then
$(\rho\circ\psi,\rho\circ\pi)$ is a Toeplitz representation of $\H$,
and since
\begin{equation*}
  (\rho\circ\pi)^{(1)}(\theta_{\xi,\eta})=(\rho\circ\psi(\xi))(\rho\circ\psi(\eta))^*=\rho\circ\pi^{(1)}(\theta_{\xi,\eta})
\end{equation*}
for all $\xi,\eta\in \H$, by linearity and continuity we have
\begin{equation*}
  (\rho\circ\pi)^{(1)}=\rho\circ\pi^{(1)}.
\end{equation*}
We denote by $\J(\H)$ the closed two-sided ideal $\phi^{-1}(\K(\H))$ in
$\A$, and we say that a Toeplitz representation $(\psi,\pi)$ of $\H$ is
\emph{Cuntz-Pimsner coinvariant} if 
\begin{equation*}
  \pi^{(1)}(\phi(a))=\pi(a)
\end{equation*}
for all $a\in \J(\H)$.

\begin{theorem}[\cite{MR1986889}*{Proposition 1.3} (cf. \cite{MR97k:46069} and
  \cite{MR2002f:46113}*{Proposition 1.6})] \label{pluss} 
  Let $\H$ be a $\cs$-correspondence over $\A$. Then there is a
  $C^*$-algebra $\Oo_\H$ and a Cuntz-Pimsner coinvariant Toeplitz
  representation $(k_\H,k_\A):\H\to \Oo_\H$ which satisfies:
  \begin{enumerate}
  \item For every Cuntz-Pimsner coinvariant Toeplitz representation
    $(\psi,\pi)$ of $\H$, there is a homomorphism $\psi\times\pi$ of
    $\Oo_\H$ such that $(\psi\times\pi)\circ k_\H=\psi$ and
    $(\psi\times\pi)\circ k_\A=\pi$,
  \item $\Oo_\H$ is generated as a $C^*$-algebra by $k_\H(\H)\cup k_\A(\A)$.
  \end{enumerate}
\end{theorem}

\begin{remark}
  The triple $(\Oo_\H,k_x,k_\A)$ is unique: if $(\X,k'_\H,k'_\A)$ has
  similar properties, then there is an isomorphism $\theta:\Oo_\H\to \X$ such
  that $\theta\circ k_\H=k'_\H$ and $\theta\circ k_\A=k'_\A$. Thus there is a
  strongly continuous gauge action $\gamma:\T\to \aut\Oo_\H$ which
  satisfies $\gamma_z(k_\A(a))=k_\A(a)$ and $\gamma_z(k_\H(x))=zk_\H(x)$
  for $a\in \A$ and $x\in \H$.
\end{remark}

\section{$\cs$-correspondences associated with subshifts}

We will now define the $\cs$-correspondence $\H_{\OSS}$ that we associate to a
 subshift $\OSS$.
%
%

We start out by defining the $C^*$-algebra $\Dt_{\OSS}$ which $\H_{\OSS}$ is a
$\cs$-correspondence over. 

\begin{definition}
For every  subshift $\OSS$ we let $\mathfrak{B}(\OSS)$ be the
abelian $C^*$-algebra of all bounded functions on $\OSS$,
and $\Dt_{\OSS}$ the $C^*$-subalgebra of
$\mathfrak{B}(\OSS)$ generated by
$\{1_{\cy\mu\nu}\mid \mu,\nu \in\a^*\}$.
\end{definition}

It turns out that the spectrum of $\Dt_{\OSS}$ is the right underlying
compact space for the $C^*$-algebra that we are going to associate
with subshifts, but since we will not need an explicit description of
this compact space we are not going to give one, but instead work with
$\Dt_{\OSS}$.


\begin{definition}
Let $\OSS$ be a  subshift. For every $a\in\a$ 
let $\Dt_a$ be the ideal in
$\Dt_{\OSS}$ generated by $1_{\tsh(C(a))}$. Let 
$\H_{\OSS}$ be the right Hilbert $\Dt_{\OSS}$-module
\begin{equation*}
\bigoplus_{a\in \a}\Dt_a
\end{equation*}
with the right action is given by 
\begin{equation*}
(f_a)_{a\in\a}f=(f_af)_{a\in\a}
\end{equation*}
and the inner product by
\begin{equation*}
\inn{(f_a)_{a\in\a},(g_a)_{a\in\a}}=\sum_{a\in\a}f_a^*g_a
\end{equation*}
for $(f_a)_{a\in\a}, (g_a)_{a\in\a}\in \bigoplus_{a\in \a}\Dt_a$ and $f\in
\Dt_{\OSS}$.   
\end{definition}

%
\begin{proposition} \label{finne}
Let $\OSS$ be a  subshift and let $a\in \a$.
Define a $*$-homomorphism
$\widetilde{\lambda}_a:\mathfrak{B}(\OSS)\to\mathfrak{B}(\OSS)$
by letting
\begin{equation*} 
  \widetilde{\lambda}_a(f)(x)=\left\{ 
    \begin{array}{ll}
      f(ax) & \textrm{if } ax\in \OSS \\
      0 & \textrm{if } ax\notin \OSS
    \end{array} \right.
\end{equation*}
for every $f\in \mathfrak{B}(\OSS)$ and every $x\in \OSS$.

Then 
$\widetilde{\lambda}_a(\Dt_{\OSS})\subseteq
\Dt_a$.
\end{proposition}

\begin{proof}
%
Let $\mu,\nu \in \a^*$ with $|\nu|\ge 1$. For every $x\in
\OSS$ is
\begin{multline*}  
    \widetilde{\lambda}_a\left(1_{\cy\mu\nu}\right)(x) = 
    \begin{cases}
      1_{\cy\mu\nu}(ax) & \text{if } ax\in \OSS \\
      0 & \text{if } ax\notin \OSS
    \end{cases}
    \\= 
    \begin{cases}
      1 & \text{if } a=\nu_1,\ x_1=\nu_2,\ \dotsc
      ,x_{|\nu|-1}=\nu_{|\nu|},\ 
      \mu\osh^{|\nu|-1}(x)\in\OSS,\ ax\in \OSS \\
      0 & \text{else.}
    \end{cases}
\end{multline*}
So $\widetilde{\lambda}_a\left(1_{\cy\mu\nu}\right)=0$
if $a\ne \nu_1$, and  
\begin{equation*}
\widetilde{\lambda}_a\left(1_{\cy\mu\nu}\right)=
1_{\cy\mu{\nu_2\nu_3\dotsm\nu_{|\nu|}}}1_{\tsh(C(a))}
\end{equation*}
if $a=\nu_1$. Hence
$\widetilde{\lambda}_a\left(1_{\cy\nu\mu}\right)\in
\widetilde{D}_a$. In a similar way, we see that
$\widetilde{\lambda}_a\left(1_{\cy\mu\emptyword}\right) =
1_{\cy{a\mu}\emptyword}$, so
$\widetilde{\lambda}_a\left(1_{\cy\nu\emptyword}\right)\in \widetilde{D}_a$.
Thus $\widetilde{\lambda}_a(\Dt_{\OSS})\subseteq \Dt_a$,
since $\widetilde{D}_{\OSS}$ is generated by
$\{1_{\cy\mu\nu}\mid \mu,\nu\in \a^*\}$.
\end{proof}

\begin{definition}
Let $\OSS$ be a subshift. We let $\phi:\Dt_{\OSS}\to \L(\H_{\OSS})$ be the $*$-homomorphism
defined by
\begin{equation*}
  \phi(f)((f_a)_{a\in \a})=(\widetilde{\lambda}_a(f)f_a)_{a\in \a}
\end{equation*}
for every $f\in \Dt_{\OSS}$ and every $(f_a)_{a\in \a}\in
\H_{\OSS}$. With this $\H_{\OSS}$ becomes a $\cs$-correspondence. 
\end{definition}

\section{The $\cs$-algebra associated with a subshift}

We are now ready to define the $\cs$-algebra $\Oo_{\OSS}$ associated with a
subshift $\OSS$.

\begin{definition}
  Let $\OSS$ be a subshift. The $\cs$-algebra $\Oo_{\OSS}$ associated
  with $\OSS$ is the $\cs$-algebra $\Oo_{\H_{\OSS}}$ from Theorem
  \ref{pluss}, where $\H_{\OSS}$ is the $\cs$-correspondence defined above. 
\end{definition}

\noindent We will now take a closer look at 
$\Oo_{\OSS}$. 
First, we show that $\Oo_{\OSS}$ is unital.

\begin{lemma} \label{unitp}
  Let $\OSS$ be a subshift and let $1=1_{\cy{\emptyword}{\emptyword}}$
    be the unit of $\widetilde{D}_{\OSS}$. Then $\kd(1)$ is a unit for
    $\Oo_{\OSS}$. 
\end{lemma}

\begin{proof}
  We have that 
  \begin{equation*}
    \kh(\xi)\kd(1)=\kh(\xi 1)=\kh(\xi),
  \end{equation*}
  and 
  \begin{equation*}
    \kd(1)\kh(\xi)=\kh(\phi(1)\xi)=\kh(\xi)
  \end{equation*}
  for every $\xi\in \H_{\OSS}$. Since we also have that
  \begin{equation*}
    \kd(1)\kd(f)=\kd(f)\kd(1)=\kd(f)
  \end{equation*} 
  for every $f\in \Dt_{\OSS}$, and $\Oo_{\OSS}$ is generated by
  $\kh(\H_{\OSS})\cup\kd(\Dt_{\OSS})$, we have that $\kd(1)$ is a unit
  for $\Oo_{\OSS}$. 
\end{proof}

\noindent We will denote the unit of $\Oo_{\OSS}$ by $I$.
 
\begin{definition}
Let $\OSS$ be a subshift. For every $a\in \a$ let $\xi_a$ be the
element $(f_{a'})_{a'\in \a}\in \H_{\OSS}$ where
$f_a=1_{\tsh(C(a))}$ and $f_{a'}=0$ for $a'\ne a$, and let for every
$\mu\in \a^*$, $S_\mu$ be the product
$\kh(\xi_{\mu_1})\kh(\xi_{\mu_2})\cdots \kh(\xi_{\mu_{|\mu|}})\in
\Oo_{\OSS}$ with the convention that $S_\emptyword=I$.
\end{definition}

\begin{lemma} \label{unit}
Let $\OSS$ be a  subshift. Let $1$ be the unit of
$\widetilde{D}_{\OSS}$, and $\Id$ the unit of $\L(\H_{\OSS}')$. Then 
\begin{equation*}
  \kd(1)=\kd^{(1)}(\Id)=\sum_{a\in \a}S_aS_a^*
\end{equation*}
is the unit of $\Oo_{\OSS}$.
\end{lemma}

\begin{proof}
It is easy to check that 
\begin{equation*}
  \phi(1)=\Id=\sum_{a\in \a}\theta_{\xi_a, \xi_a^*},
\end{equation*}
so since $(\kh,\kd)$ is a Cuntz-Pimsner coinvariant representation, 
\begin{equation*}
  \kd(1)=\kd^{(1)}(\Id)=\sum_{a\in \a}S_aS_a^*,
\end{equation*}
and we know from Lemma \ref{unitp} that $\kd(1)$ is the unit of $\Oo_{\OSS}$.
\end{proof}

\begin{lemma} \label{mayer}
Let $\OSS$ be a subshift. Then
\begin{equation*}
  \kd\left(1_{\cy\mu\nu}\right)=S_\nu S_\mu^*S_\mu S_\nu^*
\end{equation*}
for every $\mu,\nu\in \a^*$.
\end{lemma}

\begin{proof}
  Since 
  \begin{eqnarray*}
    S_a^*S_a&=&
    \kh(\xi_a)^*\kh(\xi_a)\\
    &=& \kd(\langle \xi_a,\xi_a \rangle)\\
    &=& \kd\left(1_{\tsh(C(a))}\right),
  \end{eqnarray*}
  and
  \begin{eqnarray*}
    S_a^*\kd\left(1_{\tsh^{|\mu'|}(C(\mu'))}\right)S_a &=&
    \kh(\xi_a)^*\kh\left(\phi'\left(
        1_{\tsh^{|\mu'|}(C(\mu'))}\right)\xi_a\right) \\
    &=& \kh(\xi_a)^*\kh\left(\xi_a
      \widetilde{\lambda}_a\left(1_{\tsh^{|\mu'|}(C(\mu'))}\right)\right) \\
    &=& \kh(\xi_a)^*\kh(\xi_a)\kd\left(
      \widetilde{\lambda}_a\left(1_{\tsh^{|\mu'|}(C(\mu'))}\right)\right) \\
    &=& \kd\left(1_{\tsh(C(a))}
      \widetilde{\lambda}_a\left(1_{\tsh^{|\mu'|}(C(\mu'))}\right)\right) \\
    &=& \kd\left(1_{\tsh^{|\mu' a|}(C(\mu' a))}\right)
  \end{eqnarray*}
  for every $a\in \a$ and every $\mu'\in \a^*$, we have that 
  \begin{equation*}
    S_\mu^*S_\mu=\kd\left(1_{\tsh^{|\mu|}(C(\mu))}\right)
  \end{equation*}
  for every $\mu \in \a^*$.

  It is easy to check that for every $f\in \Dt_{\OSS}$ is 
  \begin{equation*}
    \phi(f)=\sum_{a\in
      \a}\theta_{\xi_a\widetilde{\lambda}_a(f),\xi_a^*}.
  \end{equation*}
  Let $\mu,\nu\in \a^*$ with $|\nu|\ge 1$ and $a\in \a$. Then as
  proved in the proof of Proposition \ref{finne} 
  \begin{equation*}
    \widetilde{\lambda}_a\left(1_{(\cy\mu\nu}\right)=0
  \end{equation*}
  if $a\ne \nu_1$, and  
  \begin{equation*}
    \widetilde{\lambda}_a\left(1_{\cy\mu\nu}\right)=
    1_{\cy\mu{\nu_2,\dots \nu_{|\nu|}}}1_{\tsh(C(a))} 
  \end{equation*}
  if $a=\nu_1$.
  
  So 
  \begin{eqnarray*}
    \kd\left(1_{\cy\mu\nu}\right)&=&
    \kd^{(1)}\left(\phi\left(1_{\cy\mu\nu}\right)\right)\\
    &=& \kd^{(1)}\left(\theta_{\xi_{\nu_1}1_{\cy\mu{\nu_2,\dots \nu_{|\nu|}}},\xi_{\nu_1}}\right)\\
    &=& \kh(\xi_{\nu_1}1_{\cy\mu{\nu_2,\dots \nu_{|\nu|}}})\kh(\xi_{\nu_1})^*\\
    &=& S_{\nu_1}\kd(1_{\cy\mu{\nu_2,\dots \nu_{|\nu|}}})S_{\nu_1}^*.
  \end{eqnarray*}
  Hence 
  \begin{equation*}
    \kd\left(1_{\cy\mu\nu}\right)=S_\nu S_\mu^*S_\mu S_\nu^*
  \end{equation*}
  for all $\mu,\nu\in \a^*$.
\end{proof}

\begin{proposition} \label{gff}
  Let $\OSS$ be a subshift. Then $\Oo_{\OSS}$ is
  generated by $\{S_a\}_{a\in \a}$.
\end{proposition}

\begin{proof}
  $\Oo_{\OSS}$ is by Theorem \ref{pluss} generated by
  $\kh(\H_{\OSS})\cup\kd(\Dt_{\OSS})$. 

  First notice that $\kd(1)=\sum_{a\in \a}S_aS_a^*$ is in the
  $C^*$-algebra generated by $\{S_a\}_{a\in \a}$. Since 
  \begin{equation*}
    \kd\left(1_{\cy\mu\nu}\right)=S_\nu S_\mu^*S_\mu S_\nu^*
  \end{equation*}
  for all $\mu,\nu\in \a^*$, 
  and
  $\Dt_{\OSS}$ is generated by
  $\{1_{\cy\mu\nu}\mid\mu,\nu\in \a^*\}$, we have that
  $\kd(\Dt_{\OSS})$ is in the $C^*$-algebra generated by
  $\{S_a\}_{a\in \a}$. 

  Let $(f_a)_{a\in \a}\in \H_{\OSS}$. Then 
  \begin{equation*}
    (f_a)_{a\in \a}=\sum_{a\in\a}\xi_af_a,
  \end{equation*}
  so
  $\kh((f_a)_{a\in \a})=\sum_{a\in \a}S_a\kd(f_a)$, and 
  $\kh((f_a)_{a\in \a})$ is
  in the $C^*$-algebra generated by $\{S_a\}_{a\in \a}$. Hence 
  $\Oo_{\OSS}$ is
  generated by $\{S_a\}_{a\in \a}$.
\end{proof}

\section[$C^*$-algebras generated by partial isometries]{The structure
  of $C^*$-algebras generated by partial isometries} \label{bordt}

We have now established that $\Oo_{\OSS}$ 
is a unital $C^*$-algebras generated by partial isometries $\{S_a\}_{a\in
\a}$, which by Lemma \ref{unit} and \ref{mayer} satisfy

\begin{eqnarray}
\sum_{a\in \a}S_aS_a^*&=&I, \label{1} \\
S_\mu^* S_\mu S_\nu S_\nu^*&=&S_\nu S_\nu^* S_\mu^* S_\mu, \label{2} \\
S_\mu^*S_\mu S_\nu^*S_\nu&=& S_\nu^*S_\nu S_\mu^*S_\mu, \label{3}
\end{eqnarray}
where $S_\mu = S_{\mu_1}\cdots S_{\mu_{|\mu|}}$ and $S_\nu =
S_{\nu_1}\cdots S_{\nu_{|\nu|}}$, for every $\mu,\nu \in \a^*$.

We will now take a closer look at unital $C^*$-algebras generated by partial
isometries $\{S_a\}_{a\in \a}$ that satisfy the 3 relations above.

So in the rest of this section, $\a$ will be an alphabet and $\Oo$
will be a unital $C^*$-algebra generated by partial isometries
$\{S_a\}_{a\in \a}$, which satisfy the relations (\ref{1}), (\ref{2}) and
(\ref{3}) above.
  

\begin{lemma}
For every $\mu\in \a^*$, $S_\mu$ is a partial isometry.
\end{lemma}

\begin{proof}
We will prove the lemma by induction over the length of $|\mu|$. If
$|\mu|=1$, then $S_\mu$ is a partial isometry by definition. Assume
now that $S_\nu$ is a partial isometry and $a\in \a$. Then
\begin{eqnarray*}
S_{\nu a}S_{\nu a}^*S_{\nu a} &=& S_\nu S_a S_a^* S_\nu^*S_\nu S_a\\
&=& S_\nu S_\nu^*S_\nu S_a S_a^*S_a \\
&=& S_\nu S_a\\
&=& S_{\nu a}.
\end{eqnarray*}
So $S_{\nu a}$ is a partial isometry. Hence $S_\mu$ is a partial
isometry for every $\mu\in \a^*$.
\end{proof}

For $\mu\in \a^*$ we set $A_\mu=S_\mu^*S_\mu$.
We notice that since $\sum_{a\in \a}S_aS_a^*=I$, the projections
$\{S_aS_a^*\}_{a\in\a}$ are mutually orthogonal, so
$S_aS_a^*S_bS_b^*=0$ for $a\ne b$.

\begin{lemma} \label{summ}
Let $\mu,\ \nu\in \a^*$ with $|\mu|=|\nu|$. If $S_\mu^*S_\nu\ne 0$, then 
$\mu=\nu$ and $S_\mu^*S_\nu=A_\mu$.
\end{lemma}

\begin{proof}
We will prove the lemma by induction over the length of $\mu$ and $\nu$. If the length is 1 and $\mu \ne \nu$, then 
\begin{eqnarray*}
S_\mu^*S_\nu &=& S_\mu^*S_\mu S_\mu^* S_\nu S_\nu^* S_\nu \\
&=& 0
\end{eqnarray*}
since $S_\mu S_\mu^* S_\nu S_\nu^*=0$. So since $S_\mu^*S_\nu \ne 0$,
we have that $\mu=\nu$ and $S_\mu^*S_\nu=A_\mu$.

Now assume that we have proved the lemma in case $|\mu|=|\nu|=n$, and assume that $|\mu '|=|\nu '|=n+1$ and $S_{\mu '}^*S_{\nu '}\ne 0$. Set $\mu=(\mu_1',\ldots ,\mu_n')$ and $\nu=(\nu_1',\ldots ,\nu_n')$. Then $S_\mu^*S_\nu \ne 0$, so $\mu =\nu$. Since
\begin{eqnarray*}
0 &\ne& S_{\mu'}^*S_{\nu'}\\
&=& S_{\mu_{n+1}'}^*S_\mu^*S_\nu S_{\nu_{n+1}'} \\
&=& S_{\mu_{n+1}'}^*S_{\mu_{n+1}'}S_{\mu_{n+1}'}^*S_\mu^* S_\mu S_{\nu_{n+1}'}S_{\nu_{n+1}'}^*S_{\nu_{n+1}'} \\
&=& S_{\mu_{n+1}'}^*S_{\mu_{n+1}'}S_{\mu_{n+1}'}^*S_{\nu_{n+1}'}S_{\nu_{n+1}'}^*S_\mu^*S_\mu S_{\nu_{n+1}'},
\end{eqnarray*}
we have that $\mu_{n+1}'=\nu_{n+1}'$, and hence $\mu'=\nu'$. So the
lemma is true.
\end{proof}

For each $l\in \No$ we denote by $\A_l(\Oo)$ the $C^*$-subalgebra of
$\Oo$ generated by $\{A_\mu\}_{\mu \in \a_l}$. Since $\A_l(\Oo)$ is
generated by a finite number of mutually commuting projection, there
exist a finite number of mutually orthogonal projections $E_i^l,\ i=1,\ldots m(l)$, such that $(E_i^l)_{i=1, \ldots m(l)}$ is a basis for $\A_l(\Oo)$.   
We have that $\overline{\bigcup_{l\in \No}\A_l(\Oo)}$ is the $C^*$-algebra
generated by $\{A_\mu\}_{\mu\in \a^*}$. We denoted this $C^*$-algebra
by $\A(\Oo)$. Since $\A_l(\Oo)$ is finite dimensional and
$\A_l(\Oo)\subseteq \A_{l+1}(\Oo)$ for every $l\in
\No$, $\A(\Oo)$ is an AF-algebra.

\begin{lemma} \label{milan}
  For $1\le k\le l,\ \mu\in \a^k$ and $i\in \{1,2,\dots ,m(l)\}$, the
  following two conditions are equivalent:
  \begin{itemize}
  \item[a)] $S_\mu E^l_i S_\mu^* \ne 0$,
  \item[b)] $A_\mu E^l_i \ne 0$.
  \end{itemize}
\end{lemma}

\begin{proof}
  Since 
  \begin{equation*} 
    S_\mu E^l_i S_\mu^*=S_\mu A_\mu E^l_i S_\mu^*,
  \end{equation*}
  and 
  \begin{equation*}
    A_\mu E^l_i = S_\mu^* S_\mu E^l_i S_\mu^* S_\mu,
  \end{equation*}
  we have that 
  \begin{equation*}
    S_\mu E^l_i S_\mu^* \ne 0 \Leftrightarrow A_\mu E^l_i \ne 0.
  \end{equation*}
\end{proof}

\begin{lemma} \label{basis}
  Let $1\le k\le l$. Then
  \begin{itemize}
  \item[a)] For $i,i'\in \{1,2,\dots m(l)\}$ and $\mu,\mu'\in \a^k$ is
    \[ S_\mu E_i^l S_\mu^*S_{\mu'} E_{i'}^l S_{\mu'}^*= 
    \left\{ \begin{array}{ll} S_\mu E_i^l S_{\mu}^*
        &\textrm{if } \mu=\mu' \textrm{ and } i=i'\\
        0 & \textrm{if } \mu\ne \mu' \textrm{ or } i\ne i'. 
      \end{array} \right. \]
  \item[b)] $(S_\mu E_i^l S_\mu^*)^*=S_\mu E_i^l S_\mu^*$ for  $i\in
    \{1,2,\dots ,m(l)\}$ and $\mu\in\a^k$. 
  \end{itemize}
\end{lemma}

\begin{proof}
  a): By Lemma \ref{summ} 
  \begin{equation*} 
    \begin{split}
      S_\mu E_i^l S_\mu^*S_{\mu'} E_{i'}^l S_{\mu'}^*&= 
      \begin{cases}
        S_\mu E_i^l A_\mu E_{i'}^l S_{\mu'}^* &\text{if } \mu=\mu'\\
        0 & \text{if } \mu\ne \mu'
      \end{cases}\\
      &=\begin{cases}
        S_\mu A_\mu E_i^l E_{i'}^l S_{\mu'}^* &\text{if } \mu=\mu'\\
        0 & \text{if } \mu\ne \mu'
      \end{cases}\\
      &=\begin{cases}
        S_\mu E_i^l S_{\mu'}^* &\text{if } \mu=\mu'\text{ and } i=i'\\
        0 & \text{if } \mu\ne \mu'  \text{ or } i\ne i'.
      \end{cases}
    \end{split}
  \end{equation*}
        
  b): Obviously.
\end{proof}

\section{The universal property of $\Oo_{\OSS}$} \label{cole}

We let $\widetilde{\A}_{\OSS}$ be the $\cs$-subalgebra of $\Dt_{\OSS}$
generated by $\{1_{\osh^{|\mu|}(\cy{}{\mu})}\mid \mu\in\a^*\}$. 
 
\begin{lemma} \label{udvid}
Let $\OSS$ be a  subshift, $\X$ a $C^*$-algebra,
$\psi:\widetilde{\A}_{\OSS}\to \X$ a $*$-homomorphism and $\{S_a\}_{a\in
\a}$ partial isometries in $\X$ such that
\begin{itemize}
\item[a)] $\sum_{a\in \a}S_aS_a^*=\psi(1)$,
\item[b)] $S_\mu^*S_\mu S_\nu S_\nu^*=S_\nu S_\nu^* S_\mu^* S_\mu$, 
\item[c)] $S_\mu^* S_\mu=\psi\left(1_{\tsh^{|\mu|}(C(\mu))}\right)$,
\end{itemize}
where $S_\mu=S_{\mu_1}S_{\mu_2}\cdots S_{\mu_{|\mu|}}$ and
$S_\nu=S_{\nu_1}S_{\nu_2}\cdots S_{\nu_{|\nu|}}$, for every 
$\mu,\nu\in \a^*$.

Then $\psi$ extends to a $*$-homomorphism from $\Dt_{\OSS}$
to $\X$, such that 
\begin{equation*}
  \psi\left(1_{\cy\mu\nu}\right)=S_\nu S_\mu^* S_\mu S_\nu^*
\end{equation*}
for every $\mu,\nu \in \a^*$.
\end{lemma}

\begin{proof}
Let $\Oo$ be the $C^*$-subalgebra of $\X$ generated by $\{S_a\}_{a\in\a}$.
Since $S_\mu^* S_\mu=\psi\left(1_{\tsh^{|\mu|}(C(\mu))}\right)$ and
$S_\nu^* S_\nu=\psi\left(1_{\tsh^{|\nu|}(C(\nu))}\right)$, we have
that $S_\mu^*S_\mu S_\nu^*S_\nu= S_\nu^*S_\nu S_\mu^*S_\mu$ for
every $\mu,\nu \in \a^*$. Since 
\begin{eqnarray*}
S_a\psi(1)&=&S_aS_a^*S_a\psi(1)\\
&=& S_a \psi\left(1_{C(a)}\right)\psi(1)\\
&=& S_a \psi\left(1_{C(a)}\right)\\
&=& S_aS_a^*S_a\\
&=& S_a
\end{eqnarray*}
and 
\begin{eqnarray*}
\psi(1)S_a &=& \psi(1) S_aS_a^*S_a\\
&=& \sum_{a'\in \a}S_{a'}S_{a'}^*S_aS_a^*S_a \\
&=& S_aS_a^*S_a\\
&=& S_a
\end{eqnarray*}
for every $a\in \a$, $\psi(1)$ is a unit for $\Oo$. Hence $\Oo$ is
generated by partial isometries
$\{S_a\}_{a\in \a}$ which satisfy the relations (\ref{1}), (\ref{2})
and (\ref{3}) of section \ref{bordt}.

For each $l\in \No$, denote by $\widetilde{\A}_l$ the $C^*$-subalgebra
of $\widetilde{\A}_{\OSS}$ generated by $\{1_{\tsh^{|\mu|}(C(\mu))}\mid
\mu\in \a_l\}$. Since $\widetilde{\A}_l$ is generated by a
finite number of mutually commuting projections, there exists a finite
number $m(l)$ of mutually disjoint subsets $\E_i^l,\ i=1,2,\dots
,m(l)$ of $\OSS$ such that 
\begin{equation*}
  \left\{1_{\E_i^l}\mid i\in \{1,2,\dotsc ,m(l)\}\right\}
\end{equation*}
is a basis for $\widetilde{\A}_l$.

Then $\psi\left(1_{\E_i^l}\right),\ i=1,2,\dots ,m(l)$ are mutually orthogonal
projections in $\Oo$ and $\spa\left\{\psi(1_{\E_i^l})\mid i\in \{1,2,\dots
,m(l)\}\right\}=\A_l(\Oo)$. So by Lemma \ref{basis} we have that for $1\le
k\le l$ are $S_\nu \psi(1_{\E_i^l})S_\nu^*,\ i=1,2,\dots m(l)$ mutually
orthogonal projections in $\Oo$.

For each $1\le
k\le l$ denote by $\Dt_k^l$ the $C^*$-subalgebra of
$\Dt_{\OSS}$ generated by $\{1_{\cy\mu\nu}\mid \nu\in
\a^k,\mu\in \a_l\}$. It is easy to check that 
\begin{equation*}
  1_{C(\nu)\cap \tsh^{-|\nu|}(\E_i^l)},\ \nu\in \a^k,\
  i=1,2,\dots,m(l)
\end{equation*}
are mutually orthogonal projections in $\Dt_k^l$, and since 
\begin{eqnarray*}
1_{C(\nu)\cap \tsh^{-|\nu|}(\E_i^l)}=0 &\Rightarrow&
1_{\tsh^{|\nu|}(C(\nu))}1_{\E_i^l}=0\\
&\Rightarrow& S_\nu^*S_\nu \psi(1_{\E_i^l})=0 \\
&\Rightarrow& S_\nu\psi(1_{\E_i^l})S_\nu^*=0,
\end{eqnarray*}
there exists a $*$-homomorphism $\psi_k^l:\Dt_k^l\to \X$ such
that $\psi_k^l\left(1_{C(\nu)\cap
  \tsh^{-|\nu|}(\E_i^l)}\right)=S_\nu\psi(1_{\E_i^l})S_\nu^*$ for every
$\nu\in \a^k$ and every $i\in \{1,2,\dots,m(l)\}$ and hence $\psi_k^l\left(1_{\cy\mu\nu}\right)=S_\nu S_\mu^*S_\mu S_\nu^*$
for every $\nu \in \a^k$ and every $\mu \in \a_l$. 

For every $k\in \No$ denote by $\Dt_k$ the $C^*$-subalgebra
of $\Dt_{\OSS}$ generated by $\{1_{\cy\mu\nu}\mid \nu\in
\a^k,\mu\in \a^*\}$. Then
\begin{equation*}
  \Dt_k=\overline{\bigcup_{l\ge k}\Dt_k^l}.
\end{equation*}
Let $\iota_k^l$ denote the inclusion of $\Dt_k^l$ into
$\Dt_k^{l+1}$. Since $\psi_k^{l+1}\circ \iota_k^l=\psi_k^l$
for every $l\ge k$, the $\psi_k^l$'s induce a $*$-homomorphism
$\psi_k:\Dt_k \to\Oo$ such that $\psi_k\left(1_{\cy\mu\nu}\right)=S_\nu S_\mu^*S_\mu S_\nu^*$
for every $\nu \in \a^k$ and every $\mu \in \a^*$.

Since 
\begin{equation*}
  \cy\mu\nu= \bigcup_{a\in\a}\cy{\mu a}{\nu a}
\end{equation*}
for every $\mu,\nu\in \a^*$, $\Dt_k\subseteq
\Dt_{k+1}$ for every $k\in \No$ and the inclusion $\iota_k$ of
$\Dt_k$ into $\Dt_{k+1}$ is given by 
\begin{equation*}
  \iota_k\left(1_{\cy\mu\nu}\right)= \sum_{a\in \a}1_{\cy{\mu a}{\nu
      a}}.
\end{equation*}
Hence $\psi_{k+1}\circ \iota_k=\psi_k$ and since
$\Dt_{\OSS}=\overline{\bigcup_{k\in \No}\Dt_k}$,
the $\psi_k$'s induce a $*$-homomorphism $\psi:\Dt_{\OSS}\to \Oo\subseteq\X$
such that 
\begin{equation*}
  \psi\left(1_{\cy\mu\nu}\right)=S_\nu S_\mu^* S_\mu S_\nu^*
\end{equation*}
for every $\mu,\nu \in \a^*$.
\end{proof}

We are now ready to state and prove the universal property of $\Oo_{\OSS}$.

\begin{theorem} \label{unihomo}
  Let $\OSS$ be a  subshift. Then $\Oo_{\OSS}$ is the universal
  unital $\cs$-algebra generated by partial isometries $\{S_a\}_{a\in \a}$
  satisfying 
  \begin{itemize}
  \item[a)] $\sum_{a\in\a}S_aS_a^*=I$,
  \item[b)] $S_\mu^*S_\mu S_\nu S_\nu^*=S_\nu S_\nu^*S_\mu^*S_\mu$,
  \item[c)] the map $1_{\cyl{}\mu}\mapsto S_\mu^*S_\mu$ extends to a unital
    $*$-homomorphism from $\widetilde{\A}_{\OSS}$ to the $\cs$-algebra
    generated by $\{S_a\}_{a\in\a}$,
  \end{itemize}
  where $S_\mu=S_{\mu_1}\dotsm S_{\mu_{|\mu|}}$ and $S_\nu=S_{\nu_1}\dotsm
  S_{\nu_{|\nu|}}$ for every $\mu,\nu \in \a^*$.
\end{theorem}

\begin{proof}
  It follows from Lemma \ref{unit} and \ref{mayer} and Proposition
  \ref{gff} together with the
  fact that $\widetilde{\A}_{\OSS}$ is a $\cs$-subalgebra of $\Dt_{\OSS}$, that
  $\Oo_{\OSS}$ is generated by partial isometries $\{S_a\}_{a\in \a}$
  satisfying a), b) and c).
 
  Assume now that $\X$ is a unital $\cs$-algebra generated by partial
  isometries $\{T_a\}_{a\in \a}$ and that $\pi:\widetilde{\A}_{\OSS}\to
  \X$ is a unital $*$-homomorphism such that
  \begin{itemize}
  \item[a)] $\sum_{a\in \a}T_aT_a^*=I$,
  \item[b)] $T_\mu^*T_\mu T_\nu T_\nu^*=T_\nu T_\nu^* T_\mu^* T_\mu$, 
  \item[c)] $T_\mu^* T_\mu=\pi\left(1_{\tsh^{|\mu|}(\cy{}{\mu})}\right)$,
  \end{itemize}
  where $T_\mu=T_{\mu_1}T_{\mu_2}\cdots T_{\mu_{|\mu|}}$ and
  $T_\nu=T_{\nu_1}T_{\nu_2}\cdots T_{\nu_{|\nu|}}$, for every 
  $\mu,\nu\in \a^*$.

  By Lemma \ref{udvid}, $\pi$ extends to a $*$-homomorphism from $\Dt_{\OSS}$
  to $\X$, such that 
  \begin{equation*}
    \pi\left(1_{\cy\mu\nu}\right)=T_\nu T_\mu^* T_\mu T_\nu^*
  \end{equation*}
  for every $\mu,\nu \in \a^*$.
  Let
  \begin{equation*}
    \psi((f_a)_{a\in \a})= \sum_{a\in\a}T_a\pi(f_a)
  \end{equation*}
  for every $(f_a)_{a\in \a}\in \H_{\OSS}$. We will show that
  $(\psi,\pi)$ is a Cuntz-Pimsner coinvariant representation of $\H_{\OSS}$.

  It is easy to check that $\alpha \psi(\xi) +\beta
  \psi(\zeta)=\psi(\alpha \xi+\beta \zeta)$ for 
  every $\alpha,\beta \in \C$ and every $\xi,\zeta\in \H_{\OSS}$, and
  $\psi(\xi) \pi(f)=\psi(\xi f)$ for every $\xi\in
  \H_{\OSS}$ and every $f\in \widetilde{\A}_{\OSS}$. 

  Recall from the
  proof of Proposition \ref{finne} that for $\mu,\nu \in \a^*$ with
  $|\nu|\ge 1$ is $\widetilde{\lambda}_a\left(1_{\cy\mu\nu}\right)=0$
  if $a\ne \nu_1$, and $\widetilde{\lambda}_a\left(1_{\cy\mu\nu}\right) =
  1_{\cy\mu{\nu_2\nu_3\dotsm\nu_{|\nu|}}}1_{\tsh(C(a))}$ if
  $a=\nu_1$. Thus
  \begin{eqnarray*}
    \pi\left(1_{\cy\mu\nu}\right)\psi((f_a)_{a\in \a}) &=&
    \pi\left(1_{\cy\mu\nu}\right) \sum_{a\in \a}T_a\pi(f_a) \\
    &=& \sum_{a\in\a}T_\nu T_\mu^*T_\mu T_\nu^*T_a\pi(f_a) \\
    &=& T_\nu T_\mu^*T_\mu T_\nu^*T_{\nu_1}\pi(f_{\nu_1}) \\
    &=& T_{\nu_1}T_{\nu_2\nu_3\dotsm \nu_{|\nu|}} T_\mu^*T_\mu
      T_{\nu_2\nu_3\dotsm \nu_{|\nu|}}^*\pi(f_{\nu_1}) \\
    &=& T_{\nu_1}\pi\left(1_{\cy{\mu}{\nu_2\nu_3\dotsm
          \nu_{|\nu|}}}\right)\pi(f_{\nu_1}) \\
    &=& \psi\left(\xi_{\nu_1}1_{\cy{\mu}{\nu_2\nu_3\dotsm
          \nu_{|\nu|}}}f_{\nu_1}\right)\\
    &=& \psi\left((\widetilde{\lambda}_a
      (1_{\cy\mu\nu})f_a)_{a\in \a}\right)  \\
    &=& \psi\left(\phi(1_{\cy\mu\nu})(f_a)_{a\in \a}\right),
  \end{eqnarray*} 
  for $(f_a)_{a\in \a}\in \H_{\OSS}$.
  We also have that
  \begin{eqnarray*}
    \pi\left(1_{\tsh^{|\mu|}(C(\mu))}\right)\psi((f_a)_{a\in \a}) &=&
    \pi\left(1_{\tsh^{|\mu|}(C(\mu))}\right) \sum_{a\in \a}T_a\pi(f_a) \\
    &=& \sum_{a\in\a}T_\mu^*T_\mu T_a\pi(f_a) \\
    &=& \sum_{a\in\a}T_aT_{\mu a}^*T_{\mu a}\pi(f_a) \\
    &=& \sum_{a\in\a}T_a\pi\left(1_{\tsh^{|\mu a|}(C(\mu a))}\right)\pi(f_a) \\
    &=& \psi\left((1_{\tsh^{|\mu a|}(C(\mu a))}f_a)_{a\in \a}\right)\\
    &=& \psi\left((\widetilde{\lambda}_a
      (1_{\tsh^{|\mu|}(C(\mu))})f_a)_{a\in \a}\right)  \\
    &=& \psi\left(\phi(1_{\tsh^{|\mu|}(C(\mu))})(f_a)_{a\in \a}\right).
  \end{eqnarray*} 
Since $\widetilde{\D}_{\OSS}$ is generated by $1_{\cy\mu\nu},\
\mu,\nu\in\a^*$, we have that 
\begin{equation*}
  \pi(f)\psi((f_a)_{a\in \a})=\psi(\phi(f)(f_a)_{a\in \a})
\end{equation*}
for every $f\in \widetilde{\D}_{\OSS}$ and every $(f_a)_{a\in \a}\in
\H_{\OSS}$. 

Since $\sum_{a\in \a}T_aT_a^*=I$, the projections $\{T_aT_a^*\}_{a\in
  \a}$ are mutually orthogonal, so
\begin{eqnarray*}
T_a^*T_{a'} &=&
T_a^*T_aT_a^*T_{a'}T_{a'}^*T_{a'}\\
&=& 0
\end{eqnarray*}
if $a\ne a'$. Thus
\begin{eqnarray*}
  \psi((f_a)_{a\in \a})^*\psi((g_a)_{a\in \a}) &=&
  \sum_{a\in \a}\pi(f_a^*)T_a^*\sum_{a'\in \a}T_{a'}\pi(g_{a'}) \\
  &=& \sum_{a\in \a}\pi(f_a^*)T_a^*T_{a}\pi(g_{a}) \\
  &=& \sum_{a\in \a}\pi(f_a^*)\pi\left(1_{\tsh(C(a))}\right)\pi(g_{a}) \\
  &=& \pi(\langle (f_a)_{a\in \a},(g_a)_{a\in \a}\rangle)
\end{eqnarray*}
for every $(f_a)_{a\in \a},(g_a)_{a\in \a}\in \H_{\OSS}$.

Finally we see that for every $f\in \widetilde{\D}_{\OSS}$ is $\phi(f)=\sum_{a\in \a}\theta_{\xi_a \widetilde{\lambda}_a(f),\xi_a}$, so
\begin{eqnarray*}
  \pi^{(1)}(\phi(f))&=&
  \sum_{a\in\a}\psi(\xi_a\widetilde{\lambda}_a(f))\psi(\xi_a)^*\\
  &=&\sum_{a\in\a}T_a\pi(\widetilde{\lambda}_a(f))T_a^* \\
  &=& \sum_{a\in \a}\pi(f)T_aT_a^* \\
&=& \pi(f).
\end{eqnarray*}

Thus $(\psi,\pi)$ is a Cuntz-Pimsner coinvariant representation of
$\H_{\OSS}$, 
so it follows from Theorem \ref{pluss} that there exists a
$*$-homomorphism $\psi\times\pi$ from $\Oo_{\OSS}$ to $\X$ such
that $\psi\times\pi(\kh((f_a)_{a\in \a}))=\psi((f_a)_{a\in\a})$ for
every $(f_a)_{a\in \a}\in \H_{\OSS}$ and hence
\begin{equation*}
  \psi\times\pi(S_a)=\psi\times\pi(\kh(\xi_a)))=\psi(\xi_a)=T_a
\end{equation*}
for every $a\in \a$.
\end{proof}

\begin{remark} \label{remarket}
  It follows from Lemma \ref{udvid} that $\Oo_{\OSS}$ also can be
  characterized as the universal $\cs$-algebra generated by partial
  isometries $\{S_a\}_{a\in \a}$ such that the map $1_{\cy\mu\nu}\mapsto
  S_\nu S_\mu^*S_\mu S_\nu^*$ extends to a $*$-homomorphism from
  $\Dt_{\OSS}$ to the $\cs$-algebra generated by $\{S_a\}_{a\in\a}$,
  where $S_\mu=S_{\mu_1}\dotsm S_{\mu_{|\mu|}}$ and $S_\nu=S_{\nu_1}\dotsm
  S_{\nu_{|\nu|}}$ for every $\mu,\nu \in \a^*$.
\end{remark}

\begin{remark} \label{tocalg}
  Condition $b)$ can be replaced by 
  \begin{itemize}
  \item[b')] $S_\mu^* S_\mu S_\nu= S_\nu S_{\mu\nu}^*S_{\mu\nu}$,
  \end{itemize}
  because $b')$ implies that $S_\mu^*S_\mu S_\nu S_\nu^*=S_\nu
  S_{\mu\nu}^*S_{\mu\nu}S_\nu^*=S_\nu S_\nu^* S_\mu^* S_\mu S_\nu
  S_\nu^*$ and $S_\nu S_\nu^* S_\mu^* S_\mu=S_\nu (S_\mu^*S_\mu
  S_\nu)^*=S_\nu (S_\nu S_{\mu\nu}^*S_{\nu\mu})^*=S_\nu S_\nu^*S_\mu^*
  S_\mu S_\nu S_\nu^*$, and thus $S_\mu^*S_\mu S_\nu S_\nu^*=S_\nu
  S_\nu^*S_\mu^* S_\mu$, and $b)$ implies that $S_\mu^*S_\mu
  S_\nu=S_\mu^*S_\mu S_\nu S_\nu^*S_\nu=S_\nu S_\nu^*S_\mu^*S_\mu
  S_\nu=S_\nu S_{\mu\nu}^*S_{\mu\nu}$. 

  Thus if we for a two-sided subshift $\Lambda$ define $\Oo_\Lambda$
  to be $\Oo_{\OSS_\Lambda}$ where $\OSS_\Lambda$ is the one-sided
  subshift of $\Lambda$ (cf. Section \ref{notation}), then
  $\Oo_\Lambda$ has the universal property \cite{MR98h:46077}*{Theorem
    4.9} and also has the right underlying compact space
  (cf. \cite{MR2000f:46084}*{Lemma 3.1}) and thus satisfy all of the
  results of
  \cites{MR1978330,MR2000e:46087,MR98h:46077,MR2000d:46082,MR2000f:46084,MR2001e:46115,MR2001g:46147,MR1645334,MR1637788,MR1716953,MR2002h:19004}
\end{remark}

\begin{remark} \label{tmr}
  It is easy to check that for a two-sided subshift $\Lambda$, the
  partial isometries of \cite{MR2091486}*{Definition 2.1} satisfy a), b) and
  c) of Theorem \ref{unihomo} (with $\OSS=\OSS_\Lambda$, the one-sided
  subshift of $\Lambda$ (cf. Section \ref{notation})). Thus there is a
  surjective $*$-homomorphism from $\Oo_{\OSS_\Lambda}$ to
  $\Oo_\Lambda$ of \cite{MR2091486}. If $\Lambda$ satisfies the condition (I)
  of \cite{MR2091486}, then this $*$-homomorphism is an isomorphism. There
  are examples of subshifts for which the $*$-homomorphism is not
  injective (an example of this is if $\Lambda$ only consists of one point).
\end{remark}

\begin{remark} \label{remark:sergei}
  In \cite{CS} a $\cs$-algebra
  $\mathcal{D}_{\OSS}\rtimes_{\alpha,\mathcal{L}}\N$ has for every
  one-sided subshift $\OSS$ been constructed by using  Exel's crossed
  product of a $\cs$-algbra of an endomorphism. It follows from
  \cite{CS}*{Remark 9} that $\Oo_{\OSS}$ and
  $\mathcal{D}_{\OSS}\rtimes_{\alpha,\mathcal{L}}\N$ are isomorphic
  for every one-sided subshift $\OSS$. 
\end{remark}

\begin{remark} \label{remark:groupoid}
  It follows from Theorem 10 of Chapter 2 of \cite{thesis} that
  $\Oo_{\OSS}$ can be constructed as the $\cs$-algebra of a groupoid. 
\end{remark}

\section{$\Oo_{\OSS}$ is an invariant} \label{cole3}

We will now show that $\Oo_{\OSS}$ is an invariant for 
subshifts. We will do that by showing that if two subshifts
$\OSS$ and $\OSSY$ are conjugate, then $\H_{\OSS}$ and
$\H_\OSSY$ are isomorphic as $\cs$-correspondences, and it then follows
that $\Oo_{\OSS}$ and $\Oo_\OSSY$ are isomorphic.

\begin{definition} \label{hism}
Let $\X$ and $\X'$ be $C^*$-algebras, $(\H,\phi)$ a $\cs$-correspondence
over $\X$ and $(\H,\phi')$ a $\cs$-correspondence over $\X'$. If there
exist a $*$-isomorphism $\psi:\X\to \X'$ and a bijective map
$T:\H\to \H$ such that
\begin{equation*}
  \langle T\xi,\zeta \rangle=\psi(\langle \xi, T^{-1}\zeta \rangle ),
\end{equation*}
and
\begin{equation*}
  T(\phi(X)\xi)=\phi'(\psi(X))(T\xi)
\end{equation*}
for all $\xi\in \H,\ \zeta\in \H',\ X\in \X$; then we say that $(T,\psi)$ is
an $\cs$-correspondence isomorphism, $(\H,\phi)$ and $(\H',\phi')$ are
isomorphic and we write $\H\cong \H'$.
\end{definition}

It easily follows from Theorem \ref{pluss} that if $(T,\psi)$ is an
$\cs$-correspondence isomorphism from $(\H,\phi)$ to $(\H',\phi')$,
then there exists a $*$-isomorphism $T\times\psi$ from $\Oo_\H$ to
$\Oo_{\H'}$ such that $T\times\psi\circ k_\H=k_{\H'}\circ T$ and
$T\times\psi\circ\ k_\X=k_{\X'}\circ\psi$.  

\begin{lemma} \label{finne2}
Let ${\OSS}$ be a subshift.
Define a $*$-homomorphism
$\widetilde{\phi}_{\OSS}:\mathfrak{B}({\OSS})\to\mathfrak{B}({\OSS})$ by letting
\begin{equation*}
  \widetilde{\phi}_{\OSS}(f)(x)=f(\sigma(x))
\end{equation*}
for every $f\in \mathfrak{B}({\OSS})$ and every $x\in {\OSS}$.

Then $\widetilde{\phi}_{\OSS}(\widetilde{\D}_{\OSS})\subseteq
\widetilde{\D}_{\OSS}$.
\end{lemma}

\begin{proof}
Let $\mu,\nu \in \a^*$. Then
\begin{equation*}
  \sigma^{-1}(\cy\mu\nu)= \bigcup_{a\in \a}\cy\mu{a\nu},
\end{equation*}
so 
\begin{eqnarray*}
  \widetilde{\phi}_{\OSS}\left(1_{\cy\mu\nu}\right)&=&
  1_{\sigma^{-1}(\cy\mu\nu)}\\
  &=& 1_{\bigcup_{a\in \a}\cy\mu{a\nu}}\\
  &=& \sum_{a\in \a}1_{\cy\mu{a\nu}}\in \widetilde{\D}_{\OSS}.
\end{eqnarray*}
Thus, since $\widetilde{\D}_{\OSS}$ is generated by $\{1_{\cy\mu\nu}\mid \mu,\nu\in \a^*\}$ and
$\widetilde{\phi}_{\OSS}$ is a $*$-homomorphism, it follows that
$\widetilde{\phi}_{\OSS}(\widetilde{\D}_{\OSS})\subseteq
\widetilde{\D}_{\OSS}$.
\end{proof}

\begin{lemma} \label{ebbe}
Let ${\OSS}$ be a subshift. Then we have:
\begin{itemize}
\item[a)] If $\E_1,\E_2$ are subsets of ${\OSS}$ such that
  $1_{\E_1},1_{\E_2}\in \widetilde{\D}_{\OSS}$, then
  $1_{\E_1\cup\E_2}\in \widetilde{\D}_{\OSS}$.
\item[b)] If $\E$ is a subset of ${\OSS}$ such that $1_\E\in
  \widetilde{\D}_{\OSS}$, then $1_{\sigma(\E)}\in
  \widetilde{\D}_{\OSS}$.
\item[c)] If $\E$ is a subset of ${\OSS}$ such that $1_\E\in
  \widetilde{\D}_{\OSS}$, then $1_{\sigma^{-1}(\E)}\in
  \widetilde{\D}_{\OSS}$.
\end{itemize}
\end{lemma}

\begin{proof} 
  $a)$ Let $\E_1,\E_2$ be subsets of ${\OSS}$ such that
  $1_{\E_1},1_{\E_2}\in \widetilde{\D}_{\OSS}$, then
  \begin{equation*}
    1_{\E_1\cup\E_2}=1_{\E_1}+1_{\E_2}-1_{\E_1}1_{\E_2}\in
    \widetilde{\D}_{\OSS}. 
  \end{equation*}
  
  $b)$ Let $\E$ be a subset of ${\OSS}$ such that $1_\E\in
  \widetilde{\D}_{\OSS}$. Set for each $a\in \a$, 
  \begin{equation*}
    \E_a=\{x\in {\OSS}\mid ax\in \E\}.
  \end{equation*}
  It is easy to check that 
  \begin{equation*}
    \sigma(\E)=\bigcup_{a\in \a}\E_a.
  \end{equation*}
  Since $1_{\E_a}=\widetilde{\lambda}_a(1_\E)\in \widetilde{\D}_{\OSS}$
  (cf. Proposition \ref{finne}), it follows from a) that $1_{\sigma(\E)}\in
  \widetilde{\D}_{\OSS}$.

  $c)$ Let $\E$ be a subset of ${\OSS}$ such that $1_\E\in
  \widetilde{\D}_{\OSS}$. It is easy to check that
  $1_{\sigma^{-1}(\E)}=\widetilde{\phi}_{\OSS}(1_\E)$, so
  $1_{\sigma^{-1}(\E)}\in \widetilde{\D}_{\OSS}$ by Lemma \ref{finne2}.
\end{proof}

\begin{proposition} \label{isocunj}
  If two subshifts $\OSS$ and $\OSSY$ are conjugate, then
  $\Dt_{\OSS}\cong \Dt_\OSSY$ and $\H_{\OSS} \cong \H_\OSSY$.
\end{proposition}

\begin{proof}
Let $\psi:\OSS\to \OSSY$ be a conjugacy. Then we can define a
$*$-isomorphism $\Psi:\mathfrak{B}(\OSSY)\to \mathfrak{B}(\OSS)$
by setting $\Psi(f)(x)=f(\psi(x))$ for
every $f\in \mathfrak{B}(\OSSY)$ and every $x\in \OSS$.

Let $\mu\in \B(\OSSY)$. Since $C_\OSSY(\mu)$ is clopen and $\psi$ is
continuous, $\psi^{-1}(C_\OSSY(\mu))$ is clopen and hence compact. So since
$C_{\OSS}(\nu),\ \nu\in \B(\OSS)$ is a basis for the topology of $\OSS$, 
there exist a finite number of words $\mu_1,\mu_2,\dots ,\mu_r\in
\B(\OSS)$ such that
\begin{equation*}
  \psi^{-1}(C_\OSSY(\mu))=\bigcup_{k=1}^rC_{\OSS}(\mu_k).
\end{equation*}

Let $\mu,\nu\in \B(\OSSY)$ and let $\mu_1,\dots \mu_r,\nu_1,\dots
,\nu_s\in \B(\OSS)$ such that 
\begin{equation*}
  \psi^{-1}(C_\OSSY(\mu))=\bigcup_{k=1}^rC_{\OSS}(\mu_k)
\end{equation*} 
and
\begin{equation*}
  \psi^{-1}(C_\OSSY(\nu))=\bigcup_{k=1}^sC_{\OSS}(\nu_k).
\end{equation*}
Since both $\psi\circ{\osh}_{\OSS}={\osh}_{\OSSY}\circ\psi$, we have that
\begin{equation*}
  \begin{split}
    \psi^{-1}(C_\OSSY(\mu,\nu))&=
    \psi^{-1}(C_\OSSY(\nu))\cap
    \tsh_{\OSS}^{-|\nu|}(\tsh_{\OSS}^{|\mu|}(\psi^{-1}(C_\OSSY(\mu)))) \\
    &=\left(\bigcup_{k=1}^s C_{\OSS}(\nu_k)\right)\bigcap
    \left(\bigcup_{k=1}^r\tsh_{\OSS}^{-|\nu|}(\tsh_{\OSS}^{|\mu|}(C_{\OSS}(\mu_j)))\right), 
  \end{split}
\end{equation*}
so it follows from Lemma \ref{ebbe} that
\begin{eqnarray*}
\Psi\left(1_{C_\OSSY(\mu,\nu)}\right) &=&
1_{\psi^{-1}(C_\OSSY(\mu,\nu))}\in 
\Dt_{\OSS}.
\end{eqnarray*}
Hence $\Psi(\Dt_\OSSY)\subseteq
\Dt_{\OSS}$. In the same way we can prove that
$\Psi^{-1}(\Dt_{\OSS})\subseteq \Dt_\OSSY$, so $\Psi(\Dt_\OSSY)=\Dt_{\OSS}$, and thus $\Psi_{|\Dt_\OSSY}:\Dt_\OSSY\to \Dt_{\OSS}$ is a
$*$-isomorphism.

Define $T:\H_\OSSY\to \H_{\OSS}$ by
\begin{equation*}
  T(f_a)_{a\in\a_\OSSY}=
  \Bigg(\sum_{a\in\a_\OSSY}\widetilde{\lambda}_b
  \left(\Psi\left(1_{C_\OSSY(a)}\right)\right) \Psi(f_a)\Bigg)_{b\in
    \a_{\OSS}}
\end{equation*}
and $S:\H_{\OSS}\to \H_\OSSY$ by
\begin{equation*}
  S(g_b)_{b\in
    \a_{\OSS}}=\Bigg(\sum_{b\in\a_{\OSS}}\widetilde{\lambda}_a
  \left(\Psi^{-1}\left(1_{C_{\OSS}(b)}\right)\right)
  \Psi^{-1}(g_b)\Bigg)_{a\in \a_\OSSY}.
\end{equation*}

Let $a\in \a_\OSSY,\ b\in \a_{\OSS}$ and $x\in \OSSY$. If $ax\in
\OSSY$ and $(\psi^{-1}(ax))_1=b$, then 
\begin{eqnarray*}
\psi^{-1}(ax)&=&(\psi^{-1}(ax))_1\tsh(\psi^{-1}(ax))\\
&=& b\psi^{-1}(\tsh(ax))\\
&=& b\psi^{-1}(x)
\end{eqnarray*}
and thus $b\psi^{-1}(x)\in \OSS$ and $(\psi(b\psi^{-1}(x)))_1=a$.

If $b\psi^{-1}(x)\in \OSS$ and $(\psi(b\psi^{-1}(x)))_1=a$, then 
\begin{eqnarray*}
\psi(b\psi^{-1}(x))&=&
(\psi(b\psi^{-1}(x)))_1\tsh(\psi(b\psi^{-1}(x)))\\
&=& a\psi(\tsh(b\psi^{-1}(x)))\\
&=& a\psi(\psi^{-1}(x))\\
&=& ax
\end{eqnarray*}
and thus $ax\in \OSSY$ and $(\psi^{-1}(ax))_1=b$.

Hence $(ax\in \OSSY \land (\psi^{-1}(ax))_1=b) \Leftrightarrow
(b\psi^{-1}(x)\in \OSS \land (\psi(b\psi^{-1}(x)))_1=a)$. So
\begin{eqnarray*}
\Psi^{-1}\left(\widetilde{\lambda}_b\left(\Psi\left(1_{C_\OSSY(a)}\right)\right)\right)(x)&=& \left\{ \begin{array}{ll}
\Psi\left(1_{C_\OSSY(a)}\right)(b\psi^{-1}(x)) & \textrm{if } b\psi^{-1}(x)\in \OSS \\
0 & \textrm{if } b\psi^{-1}(x)\notin \OSS
\end{array} \right.\\
&=&\left\{ \begin{array}{ll}
1 & \textrm{if }  b\psi^{-1}(x)\in \OSS \land (\psi(b\psi^{-1}(x)))_1=a\\
0 & \textrm{else}
\end{array} \right.\\
&=&\left\{ \begin{array}{ll}
1 & \textrm{if } ax\in \OSSY\land (\psi^{-1}(ax))_1=b \\
0 & \textrm{else}
\end{array} \right.\\
&=& \left\{ \begin{array}{ll}
\Psi^{-1}\left(1_{C_{\OSS}(b)}\right)(ax) & \textrm{if } ax\in \OSSY \\
0 & \textrm{if } ax\notin \OSSY
\end{array} \right.\\
&=& \widetilde{\lambda}_a\left(\Psi^{-1}\left(1_{C_{\OSS}(b)}\right)\right)(x)
\end{eqnarray*}
and hence
$\Psi^{-1}\left(\widetilde{\lambda}_b\left(\Psi\left(1_{C_\OSSY(a)}\right)\right)\right)=\widetilde{\lambda}_a\left(\Psi^{-1}\left(1_{C_{\OSS}(b)}\right)\right)$
for all $a\in \a_\OSSY$ and $b\in \a_{\OSS}$, and thus
\begin{eqnarray*}
\langle T(f_a)_{a\in \a_\OSSY}, (g_b)_{b\in \a_{\OSS}} \rangle &=&
\left\langle \Bigg(\sum_{a\in\a_\OSSY}\widetilde{\lambda}_b\left(\Psi\left(1_{C_\OSSY(a)}\right)\right)\Psi(f_a)\Bigg)_{b\in \a_{\OSS}},(g_b)_{b\in \a_{\OSS}} \right\rangle \\
&=& \sum_{b\in
  \a_{\OSS}}\sum_{a\in\a_\OSSY}\widetilde{\lambda}_b\left(\Psi\left(1_{C_\OSSY(a)}\right)\right)\Psi(f_a^*)g_b\\
&=&  \sum_{a\in\a_\OSSY}\Psi(f_a^*)\sum_{b\in
  \a_{\OSS}}\widetilde{\lambda}_b\left(\Psi\left(1_{C_\OSSY(a)}\right)\right)g_b\\
&=& \Psi\Bigg(\left\langle (f_a)_{a\in \a_\OSSY},
  \Bigg(\sum_{b\in\a_{\OSS}}\widetilde{\lambda}_a\left(\Psi^{-1}\left(1_{C_{\OSS}(b)}\right)\right)\Psi^{-1}(g_b)\Bigg)_{a\in \a_\OSSY}\right\rangle \Bigg) \\
&=& \Psi(\left\langle (f_a)_{a\in \a_\OSSY}, S(g_b)_{b\in
  \a_{\OSS}}\right\rangle ) 
\end{eqnarray*}
for all $(f_a)_{a\in \a_\OSSY}\in \H_\OSSY$ and all $(g_b)_{b\in
  \a_{\OSS}}\in \H_{\OSS}$.

Let $a\in \a_\OSSY,\ b\in \a_{\OSS}$ and $y\in \OSS$. If $by\in
\OSS$ and $(\psi(by))_1=a$, then
\begin{eqnarray*}
\psi(by)&=& (\psi(by))_1 \tsh(\psi(by))\\
&=& a\psi(\tsh(by))\\
&=& a\psi(y),
\end{eqnarray*}
and thus $a\psi(y)\in \OSSY$.
 
So for every $f\in \Dt_\OSSY$ is
\begin{eqnarray*} 
\widetilde{\lambda}_b\left(\Psi\left(1_{C_\OSSY(a)}\right)\right)\Psi(\widetilde{\lambda}_a(f))(y)&=&
\left\{\begin{array}{ll} \Psi\left(1_{C_\OSSY(a)}\right)(by)\Psi(\widetilde{\lambda}_a(f))(y) &
    \textrm {if } by\in \OSS \\
0 & \textrm{if } by\notin \OSS
\end{array} \right. \\
&=& \left\{\begin{array}{ll} 1_{C_\OSSY(a)}(\psi(by))\widetilde{\lambda}_a(f)(\psi(y)) &
    \textrm {if } by\in \OSS \\
0 & \textrm{if } by\notin \OSS
\end{array} \right. \\
&=& \left\{\begin{array}{ll} f(a\psi(y)) &
    \textrm {if } by\in \OSS,\ (\psi(by))_1=a \textrm{ and }
    a\psi(y)\in \OSSY \\
0 & \textrm{else}
\end{array} \right. \\
&=& \left\{\begin{array}{ll} f(\psi(by)) &
    \textrm {if } by\in \OSS \textrm{ and } (\psi(by))_1=a\\
0 & \textrm{else}
\end{array} \right. \\
&=& \widetilde{\lambda}_b\left(\Psi\left(1_{C_\OSSY(a)}\right)\right)\widetilde{\lambda}_b(\Psi(f))(y),
\end{eqnarray*}
and hence $\widetilde{\lambda}_b\left(\Psi\left(1_{C_\OSSY(a)}\right)\right)\Psi(\widetilde{\lambda}_a(f))=\widetilde{\lambda}_b\left(\Psi\left(1_{C_\OSSY(a)}\right)\right)\widetilde{\lambda}_b(\Psi(f))$
for all $f\in \Dt_\OSSY$, $a\in \a_\OSSY$ and $b\in
\a_{\OSS}$. Thus
\begin{eqnarray*}
T(\phi'(f)(f_a)_{a\in \a_\OSSY})&=& T(\widetilde{\lambda}_a(f)f_a)_{a\in
  \a_\OSSY}\\
&=& \Bigg(\sum_{a\in
    \a_\OSSY}\widetilde{\lambda}_b\left(\Psi\left(1_{C_\OSSY(a)}\right)\right)\Psi(\widetilde{\lambda}_a(f)f_a)\Bigg)_{b\in
    \a_{\OSS}}\\
&=& \Bigg(\widetilde{\lambda}_b(\Psi(f))\sum_{a\in
    \a_\OSSY}\widetilde{\lambda}_b\left(\Psi\left(1_{C_\OSSY(a)}\right)\right)\Psi(f_a)\Bigg)_{b\in \a_{\OSS}}\\
&=& \phi'(\Psi(f))T(f_a)_{a\in \a_\OSSY}
\end{eqnarray*}
for all $(f_a)_{a\in \a_\OSSY}\in \H_\OSSY$ and all $f\in
\Dt_\OSSY$.

Since $\Psi^{-1}\left(\widetilde{\lambda}_b\left(\Psi\left(1_{C_\OSSY(a)}\right)\right)\right)=\widetilde{\lambda}_a\left(\Psi^{-1}\left(1_{C_{\OSS}(b)}\right)\right)$
for all $a\in \a_\OSSY$ and $b\in \a_{\OSS}$, $\Psi,\ \Psi^{-1}$ and 
$\widetilde{\lambda}_a$ are $*$-homomorphisms and 
\begin{equation*}
  1_{C_{\OSS}(b)}1_{C_{\OSS}(b')}=
  \left\{\begin{array}{ll} 
      1_{C_{\OSS}(b)} & \textrm{if }b=b'\\
      0 & \textrm{if }b\ne b'
    \end{array} \right.
\end{equation*}
for $b,b'\in \a_{\OSS}$, we have that
\begin{equation*}
  \widetilde{\lambda}_b \left(\Psi\left(1_{C_\OSSY(a)}\right)\right)
  \widetilde{\lambda}_{b'}
  \left(\Psi\left(1_{C_\OSSY(a)}\right)\right)= \left\{
    \begin{array}{ll} 
      \widetilde{\lambda}_b\left(\Psi\left(1_{C_\OSSY(a)}\right)\right) &\textrm{if } b=b'\\ 
      0 & \textrm{if }b\ne b'
    \end{array}\right. 
\end{equation*}
for all $a\in \a_\OSSY$ and all $b,b'\in \a_{\OSS}$; and hence
\begin{eqnarray*}
TS(g_b)_{b\in \a_{\OSS}}&=&
T\Bigg(\sum_{b\in\a_{\OSS}}\widetilde{\lambda}_a\left(\Psi^{-1}\left(1_{C_{\OSS}(b)}\right)\right)\Psi^{-1}(g_b)\Bigg)_{a\in
  \a_\OSSY}\\
&=& \Bigg(\sum_{a\in
  \a_\OSSY}\widetilde{\lambda}_{b'}\left(\Psi\left(1_{C_\OSSY(a)}\right)\right)\sum_{b\in\a_{\OSS}}\Psi\left(\widetilde{\lambda}_a\left(\Psi^{-1}\left(1_{C_\OSSY(b)}\right)\right)\right)g_b\Bigg)_{b'\in \a_{\OSS}}\\
&=&\Bigg(\sum_{a\in\a_\OSSY}\widetilde{\lambda}_{b'}\left(\Psi\left(1_{C_{\OSS}(a)}\right)\right)\sum_{b\in\a_{\OSS}}\widetilde{\lambda}_b\left(\Psi\left(1_{C_{\OSS}(a)}\right)\right)g_b\Bigg)_{b'\in
  \a_{\OSS}} \\
&=&\Bigg(\sum_{a\in\a_\OSSY}\widetilde{\lambda}_b\left(\Psi\left(1_{C_\OSSY(a)}\right)\right)g_b\Bigg)_{b\in
  \a_{\OSS}} \\
&=&\Bigg(\widetilde{\lambda}_b\Bigg(\Psi\Bigg(\sum_{a\in\a_\OSSY}1_{C_\OSSY(a)}\Bigg)\Bigg)g_b\Bigg)_{b\in
  \a_{\OSS}} \\
&=& (\widetilde{\lambda}_b(1)g_b)_{b\in \a_{\OSS}}\\
&=& \left(1_{\tsh(C_{\OSS}(b))}g_b\right)_{b\in \a_{\OSS}}\\
&=& (g_b)_{b\in \a_{\OSS}}
\end{eqnarray*}
for all $(g_b)_{b\in \a_{\OSS}}\in \H_{\OSS}$.

In the same way one can prove that $ST(f_a)_{a\in
  \a_\OSSY}=(f_a)_{a\in \a_\OSSY}$ for all $(f_a)_{a\in
  \a_\OSSY}\in \H_\OSSY$.

Hence $(T,\Psi)$ is a $\cs$-correspondence isomorphism and $\H_{\OSS}
\cong \H_\OSSY$.
\end{proof} 

\begin{remark}
  With Proposition \ref{isocunj} in hand, it is not difficult to prove
  that if $\Lambda$ and $\Gamma$ are two two-sided subshifts such that
  the one-sided subshift $\OSS_\Lambda$ of $\Lambda$ and the one-sided
  subshift $\OSS_\Gamma$ of $\Gamma$ are conjugate, then $\Oo_\Lambda$
  and $\Oo_\Gamma$ of \cite{MR2091486} (cf. Remark \ref{tmr}) are isomorphic.
\end{remark}

\begin{theorem} \label{toke} 
If two subshifts $\OSS$ and $\OSSY$ are conjugate, then there exists a
$*$-isomorphism $\rho$ from $\Oo_{\OSS}$ to $\Oo_{\OSSY}$ such that
$\gamma_z\circ\rho=\rho\circ\gamma_z$ for every $z\in \T$.
\end{theorem}

\begin{proof}
It follows from Theorem \ref{isocunj} that there exists a
$\cs$-correspondence isomorphism $(T,\Psi)$ from $\H_{\OSS}$ to
$\H_{\OSSY}$. Thus there exists a $*$-isomorphism $\rho:\Oo_{\OSS}\to
\Oo_\OSSY$ such that $\rho(k_{\H_{\OSS}}(\xi))=k_{\H_{\OSSY}}(T\xi)$
for every $\xi\in \H_{\OSS}$ and
$\rho(k_{\Dt_{\OSS}}(f))=k_{\Dt_{\OSSY}}(\Psi(f))$ for every $f\in \Dt_{\OSS}$. Hence 
\begin{eqnarray*}
\gamma_z(\rho(k_{\H_{\OSS}}(\xi)))&=&\gamma_z(k_{\H_{\OSSY}}(T\xi))\\
&=& zk_{\H_{\OSSY}}(T\xi)\\
&=& k_{\H_{\OSSY}}(Tz\xi)\\
&=& \rho(k_{\H_{\OSS}}(z\xi))\\
&=& \rho(zk_{\H_{\OSS}}(\xi))\\
&=& \rho(\gamma_z(k_{\H_{\OSS}}(\xi)))
\end{eqnarray*}
for every $\xi\in \H_{\OSS}$ and every $z\in \T$, and
\begin{eqnarray*}
\gamma_z(\rho(k_{\Dt_{\OSS}}(f)))&=&\gamma_z(k_{\Dt_{\OSSY}}(\Psi(f)))\\
&=& zk_{\Dt_{\OSSY}}(\Psi(f))\\
&=& k_{\Dt_{\OSSY}}(\Psi(zf))\\
&=& \rho(k_{\Dt_{\OSS}}(zf))\\
&=& \rho(zk_{\Dt_{\OSS}}(f))\\
&=& \rho(\gamma_z(k_{\Dt_{\OSS}}(f)))
\end{eqnarray*}
for every $f\in \Dt_{\OSS}$ and every $z\in \T$. Since
$\Oo_{\OSS}$ is generated by $k_{\H_{\OSS}}(\H_{\OSS})\cup k_{\Dt_{\OSS}}(\Dt_{\OSS})$, it follows 
that $\gamma_z\circ \rho=\rho\circ \gamma_z$ for every $z\in \T$.
\end{proof}

\begin{remark}
  One can prove that if $\Lambda$ and $\Gamma$ are two two-sided
  subshifts which are flow equivalent (cf. \cite{MR53:9179},
  \cite{MR86j:58078}, \cite{MR98a:46082} and
  \cite{MR97a:58050}*{§13.6}), then $\Oo_{\OSS_\Lambda}$ and
  $\Oo_{\OSS_\Gamma}$, where $\OSS_\Lambda$ is the one-sided subshift
  of $\Lambda$ and $\OSS_\Gamma$ is the one-sided subshift of
  $\Gamma$, are stably isomorphic. This has been proved in
  \cite{MR2002h:19004} (cf. Remark \ref{tocalg}) under the assumption
  of condition (I) and (E) (cf. \cite{MR2001e:46115}), but there is an
  alternative proof, which will appear in \cite{tmcx}, which does
  not require condition (I) and (E).

  Notice that if the two two-sided subshifts $\Lambda$ and $\Gamma$
  are conjugate (as two-sided subshifts,
  cf. \cite{MR97a:58050}*{Definition 1.5.9.}), then $\Lambda$ and
  $\Gamma$ are flow equivalent (cf. \cite{MR97a:58050}*{§13.6}), and
  so $\Oo_{\OSS_\Lambda}$ and $\Oo_{\OSS_\Gamma}$ are stably
  isomorphic. One can in fact show that if $\Lambda$ and $\Gamma$
  are conjugate, then there exists a $*$-isomorphism from $\Phi$ from
  $\Oo_{\OSS_\Lambda}\otimes \K$ to  $\Oo_{\OSS_\Gamma}\otimes \K$ such
  that $\Phi(\D_{\OSS_\Lambda}\otimes \CC)=\D_{\OSS_\Gamma}\otimes
  \CC$ and such that $\Phi\circ\gamma_\Gamma\circ\Phi^{-1}$ and
  $\gamma_\Lambda$ are exterior equivalent, where $\K$ is the
  $C^*$-algebra of all compact operators on a separable
  infinite-dimensional Hilbert space, $\CC$ is a maximal commutative
  $C^*$-subalgebra of $\K$, $\D_{\OSS_\Lambda}$ (respectively
  $\D_{\OSS_\Gamma}$) is the $C^*$-subalgebra of $\Oo_\Lambda$
  (respectively $\Oo_\Gamma$) generated by $\{S_\nu S_\mu^* S_\mu
  S_\nu^*\}_{\mu,\nu\in \B(\OSS_\Lambda)}$ (respectively $\{S_\nu S_\mu^* S_\mu
  S_\nu^*\}_{\mu,\nu\in \B(\OSS_\Gamma)}$) (noticed that
  $\D_{\OSS_\Lambda}$ (respectively $\D_{\OSS_\Gamma}$) is isomorphic
  to $\Dt_{\OSS_\Lambda}$ (respectively $\Dt_{\OSS_\Gamma}$),
  cf. Remark \ref{remarket}) and $\gamma_\Lambda$ (respectively
  $\gamma_\Gamma$) is the gauge action on $\Oo_{\OSS_\Lambda}$
  (respectively $\Oo_{\OSS_\Gamma}$). This has been proved in
  \cite{MR2001e:46115}*{Corollary 6.2} for two-sided subshifts
  satisfying the conditions (I) and (E), and will be proved in full
  generality in \cite{tmcx}.

\end{remark}

\end{document}